\documentclass[12pt]{extarticle}
\usepackage{amsmath, amssymb}
\usepackage[cp1251]{inputenc}
\usepackage[russian]{babel}

\begin{document}

\begin{center}
\textbf{Steiner symmetrization and the initial coefficients of
univalent functions}\vskip0.3cm

V. N. Dubinin
\end{center}

\textbf{Abstract.} We establish the inequality $|a_{1}|^{2}-Re \
a_{1}a_{-1}\geq |a^{*}_{1} |^{2}-Re \ a^{*}_{1} a^{*}_{-1}$ for
the initial coefficients of any function $f(z) =
a_{1}z+a_{0}+a_{-1}/z+ ...$ meromorphic and univalent in the
domain $D = \{z : |z| > 1\}$, where $a^{*}_{1}$ and $a^{*}_{-1}$
are the corresponding coefficients in the expansion of the
function $f^{*}(z)$ that maps the domain $D$ conformally and
univalently onto the exterior of the result of the Steiner
symmetrization with respect to the real axis of the complement of
the set $f(D)$. The P\'{o}lya–Szeg\H{o} inequality $|a_{1}| \geq
|a^{*}_{1}|$ is already known. We describe some applications of
our inequality to functions of class $\Sigma$.

\textbf{Keywords:} Steiner symmetrization, capacity of a set,
univalent function, covering theorem.

\begin{center}
\textbf{§ 1. Introduction and statement of results}
\end{center}

For arbitrary real $u$ we denote by $l(u)$ the vertical line $Re \
w = u$. By the \textit{Steiner symmetrization} of a closed bounded
set $E \subset \mathbb{C}_{w}$ with respect to the real axis we
mean the result of transforming this set into a symmetric set
$$
E^{*} = \{w = u + iv : E \cap l(u) \neq \varnothing, \, 2|v| \leq
\mu(E \cap l(u))\},
$$
where $\mu$ stands for linear Lebesgue measure. Let $E$ be a
non-degenerate continuum and let a function
$$
f(z) = a_{1}z + a_{0} + \frac{a_{-1}}{z} + ...
$$
map the exterior $D:= \{z : |z|>1\}$ of the unit disc conformally
and univalently onto the connected component of
$\mathbb{\overline{C}}_{w} \backslash E$ containing the point at
infinity. We denote by
$$
f^{*}(z) = a^{*}_{1} z + a^{*}_{0} + \frac{a^{*}_{-1}}{z} + ...
$$
a function taking the domain $D$ conformally and univalently onto
$\mathbb{\overline{C}}_{w} \backslash E^{*}$. The following
inequality of P\'{o}lya–Szeg\H{o} [1] is well known:
$$
|a_{1}| \geq |a^{*}_{1} |,\eqno (1)
$$
and is equivalent to the following inequality for the logarithmic
capacities:
$$
\mathrm{cap}\,E\geqslant\mathrm{cap}\,E^*.
$$
There naturally arises the question of the behaviour of the
subsequent coefficients in the decomposition of $f$ under Steiner
symmetrization. There is no meaningful inequality between $|a_0|$
and $|a_0^*|$. Indeed, translation of the set $E$ along the
imaginary axis preserves the values
$\mathrm{Re}\,a_0,\;\mathrm{Re}\,a_0^*$, and
$\mathrm{Im}\,a_0^*=0$ and makes the value $\mathrm{Im}\,a_0$
quite arbitrary. On the other hand, there are examples for which
$\mathrm{Re}\,a_0\neq \mathrm{Re}\,a_0^*$, and then translation of
$E$ along the real axis leads to the inequality
$|\mathrm{Re}\,a_0|<|\mathrm{Re}\,a_0^*|$. Further, if $E$ is a
line segment forming an acute angle with the real axis, then
$|a_{-1}|>|a_{-1}^*|$. In the case when $E=\{w:|w|=1\}$, we have
the reverse inequality, $|a_{-1}|=0<1/2=|a_{-1}^*|$. We shall
prove the following result.

\textbf{Theorem 1.}\textit{The functions $f$ and $f^*$ defined
above satisfy the inequality}\vspace{0.2cm}
$$
|a_1|^2-\mathrm{Re}\,a_1a_{-1}\geqslant
|a_1^*|^2-\mathrm{Re}\,a_1^*a_{-1}^*.\eqno (2)
$$

This inequality has the following interpretation in terms of
capacity. Let a set $E$ be symmetric with respect to the real axis
and let $H\setminus E$ be a simply connected domain, where
$H:=\{w:\mathrm{Im}\,w>0\}$. We denote by $g$ the function which
takes the domain $H\setminus E$ conformally and univalently onto
the half-plane $H$ in such a way that
$$
\lim\limits_{w\rightarrow\infty}[g(w)-w]=0.
$$
The limit
$$
\mathrm{hcap}(E\cap H)=\lim\limits_{w\rightarrow\infty} w[g(w)-w]
$$
is referred to as the \textit{half-plane capacity (from infinity)}
of the set $E\cap H$ ([2], p. 69). In view of the expansion of
$f$, we conclude that $g$ has the following form in a
neighbourhood of the point at infinity:
$$
g(w)=w+\frac{a_1^2-a_1a_{-1}}{w}+\ldots\;,
$$
where $a_1$ and $a_{-1}$ are real numbers. Thus,
$$
\mathrm{hcap}(E\cap H)=|a_1|^2-\mathrm{Re}\,a_1a_{-1}
$$
and the inequality (2) can be written in the form
$$
\mathrm{hcap}(E\cap H)\geqslant \mathrm{hcap}(E^*\cap H).
$$

Theorem 1 is supplemented by the following assertion.

\vskip0.3cm \textbf{Theorem 2.} \textit{Let a function $
\widetilde{f}(z)=\widetilde{a}_1z+\widetilde{a}_0+\widetilde{a}_{-1}/z+\ldots$
map $D$ conformally and univalently onto the exterior of a
continuum $ \widetilde{E}\subset E^*$. Then}
$$
|a_1^*|^2-\mathrm{Re}\,a_1^*a_{-1}^*\geqslant
|\widetilde{a}_1|^2-\mathrm{Re}\,\widetilde{a}_1\widetilde{a}_{-1}.\eqno(3)
$$
\vskip0.3cm Proofs of Theorems 1 and 2 are given in the concluding
part of the paper. To obtain the inequality (2), we need the
symmetrization with respect to a circle [3], as described in § 2.
The inequality (3) follows from a result of Schiffer which was
established using Hadamard’s formula for the variation of the
Green function ([4], § 3). As applications of Theorems 1 and 2, we
prove covering results for the well-known class $\Sigma$ of
functions $f(z)=z+a_0+a_{-1}/z+\ldots$ that are meromorphic and
univalent in $D$ [5].

\vskip0.3cm \textbf{Corollary 1.} \textit{Let $f$ be a function
belonging to the class $\Sigma$ and let $w_0$ be an arbitrary
point of the complement $E=\mathbb{C}_w\setminus f(D)$. Then the
inequality}
$$
\frac{m_f^4(w_0,\varphi)+16R_f^4(w_0)}{
8m_f^2(w_0,\varphi)}\leqslant 1+\mathrm{Re}\,e^{-2i\varphi}a_{-1},
$$
\textit{holds for any real number $\varphi$, where
$R_f(w_0)\geqslant 0$ stands for the radius of the largest disc
centred at the point $w_0$ and belonging to the set $E$, and
$m_f(w_0,\varphi)$  is the linear Lebesgue measure of the
intersection of $E$ with the line $\{w=w_0+te^{i\varphi}:t\in
\mathbb{R}\}$. This inequality becomes equation for the functions
$f(z)=w_0+e^{i\varphi}\lambda^{-1}h^{-1}(\lambda
h(e^{-i\varphi}z))$ with $h(\zeta)=\zeta+1/\zeta$ and any
$\lambda>1$.}

In particular, the following inequalities hold:
$$
\frac{1}{8}m_f^2\left
(w_0,\frac{\mathrm{arg}\,a_{-1}}{2}\right)-1\leqslant
|a_{-1}|\leqslant
1-\frac{1}{8}m_f^2\left(w_0,\frac{\mathrm{arg}\,a_{-1}+
\pi}{2}\right).
$$
The right-hand inequality refines a well-known corollary to the
area theorem: $|a_{-1}|\leqslant 1$ ([5], Ch. II, \S4). Both
inequalities supplement the classical bound
$$
m_f(w_0,\varphi)\leqslant 4 \,\,\,\, \forall\varphi,
$$
which follows from (1). Namely,
$$
m_f \left(w_0,\frac{\mathrm{arg}\,a_{-1}}{2}\right)\leqslant
\sqrt{8(1+|a_{-1}|)}\leqslant 4,
$$
$$
m_f \left(w_0,\frac{\mathrm{arg}\,a_{-1}+ \pi}{2}\right)\leqslant
\sqrt{8(1-|a_{-1}|)}.
$$

These inequalities become equalities when $|a_{-1}|=1$ and
$f(z)=z+w_0+e^{2i\varphi}/z$. It would be of interest to obtain
sharp estimates for a fixed $|a_{-1}| \neq 1$.

\vskip0.3cm \textbf{Corollary 2.} \textit{Suppose that a function
$f(z)=z+a_0+a_{-1}/z+\ldots$ of class $\Sigma$ satisfies the
inequality}
$$
\mu\left(( \mathbb{C}_w\setminus f(D))\cap l(u)\right)\geqslant
\alpha\qquad \forall u,\qquad\beta\leqslant u\leqslant \gamma.
$$
\textit{for some $\alpha$, $\beta$ and $\gamma$. Then}
$$
\mathrm{Re}\,a_{-1}\leqslant 1-\frac{c^2}{2}(1-k^2),
$$
\textit{where the real constants $c$ and $k$ can be found from the
condition}
$$
c\int\limits_0^1\sqrt{\frac{\zeta^2-k^2}{\zeta^2-1}}d\zeta=\frac{\gamma-\beta}{2}-\frac{i\alpha}{2},\,\,\,c>0,\,\,0<k<1.
$$
\textit{This inequality becomes an equality for a function $f$ of class $\Sigma$ mapping $D$ conformally
and univalently onto the exterior of a rectangle with sides lying on the lines
$u=\beta,\;u=\gamma,\;\gamma-\beta< 4$, and of an appropriate height $\alpha$.}

Corollaries 1 and 2 are obtained by successively applying the
inequalities (2) and (3). The list of assertions of this kind can
readily be extended in the same way as the well-known applications
of Steiner symmetrization to function theory ([3],[6]).

\begin{center}\bf \S 2. Symmetrization with respect to a circle\end{center}

Following [3], § 3, we denote by $r(w)$ the regular branch of the
function $\zeta=i\log\,w$ mapping the plane $ \mathbb{C}_w$ with a
cut along the real negative semi-axis onto the strip
$-\pi<\mathrm{Re}\,\zeta<\pi$. We define the values of the
function $r(w)$ on the cut in the sense of the boundary
correspondence. Let $E$ be an arbitrary closed set in the plane
$\overline{\mathbb{C}}_w$ that does not contain both the origin
and the point at infinity. By the \textit{symmetrization of $E$
with respect to the circle} $|w|=1$ we mean the passage from $E$
to the symmetric set with respect to $|w|=1$,
$$
RE=r^{-1}((r(E))^*),
$$
where, as above, the symbol $^*$ stands for the Steiner
symmetrization with respect to the real axis, carried out in the
strip $-\pi\leqslant\mathrm{Re}\,\zeta\leqslant\pi$. We now give a
direct definition of this transformation. For a closed set $E$
lying in $\mathbb{C}_w\setminus\{0\}$, we write
$$
E(\theta)=E\cap\{w:\mathrm{arg}\,w=\theta\},\quad
R(\theta)=\exp\left(
\frac12\int_{E(\theta)}\frac{d\rho}{\rho}\right),
$$
$$
\widetilde{E}(\theta)=\left\{\begin{array}{ll} \{w=\rho
e^{i\theta}:R^{-1}(\theta)\leqslant \rho\leqslant
R(\theta)\},&\mbox{ if }E(\theta)\neq \varnothing,\\
\varnothing,&\mbox{ if } E(\theta)=\varnothing.\end{array}\right.
$$
It can readily be seen that
$$
RE=\bigcup\limits_{0\leqslant\theta\leqslant 2\pi}
\widetilde{E}(\theta).
$$
For a fixed $v>0$ we denote by $L_v(w):=w/{v-i}$ the linear
transformation taking the circle $|w-iv|=v$ to the unit circle
$|w|=1$ and $L_v^{-1}$ the inverse map. By the \textit{result of
the symmetrization} of a closed bounded set $E$ with respect to
the circle $|w-iv|=v$ we mean the set
$$
R_vE=L_v^{-1}(RL_v(E)).
$$
For an open set $B$ containing the points $iv$ and $\infty$, we
write
$$
S_vB=\overline{\mathbb{C}}_w\setminus
R_v(\overline{\mathbb{C}}_w\setminus B).
$$

\vskip0.3cm \textbf{Lemma 1.} \textit{If open sets $B_1$ and $B_2$
satisfy the conditions $\infty \in B_1$ and $\overline{B}_1\subset
B_2$, then the inclusion relation}
$$
S_vB_1\subset
\overline{\mathbb{C}}_w\setminus(\overline{\mathbb{C}}_w\setminus
B_2)^*\qquad\qquad\qquad (R_v( \overline{\mathbb{C}}_w\setminus
B_1)\supset ( \overline{\mathbb{C}}_w\setminus B_2)^*).
$$
\textit{holds for all sufficiently large} $v>0$.

\vskip0.3cm The proof of Lemma 1 is clearly of a technical nature, and therefore we omit it.
We only note the importance of the condition that $B_1$ is contained in a compact
subset of $B_2$. Then the closed set $\overline{\mathbb{C}}_w\setminus B_2$ is contained in $\overline{\mathbb{C}}_w\setminus B_1$ together
with some neighbourhood $U$. Near the real axis, the rays passing through the
point $iv$ and intersecting the neighbourhood $U$ tend to lines parallel to the imaginary
axis as $v\rightarrow\infty$. Here the ‘logarithmic measure’ in a neighbourhood of the
circle $|w-iv|=v$ tends to the Euclidean measure.

Let $g_B(z,z_0)$ denote the Green function of the connected component of $B$ which
contains the point $z_0$ (with a pole at this point), where $g_B(z,z_0)$ is defined to be
zero outside this connected component. Let $r(B,z_0)$ be the inner radius of the
above component with respect to the point $z_0$ [3].

\vskip0.3cm \textbf{Lemma 2.} \textit{If the connected components
of the open set $B$ have Green’s functions and the points $iv$ and
$\infty$ belong to $B$, then}
$$
\log[r(B,iv)r(B,\infty)]+2g_B(iv,\infty)\leqslant
\log[r(S_vB,iv)r(S_vB,\infty)]+2g_{S_vB}(iv,\infty)
$$
\textit{for any $v>0$.} \textit{Proof.} This follows from [3],
Theorem 1.7, Proposition 1.11 (see also [7], Theorem 1).

\begin{center}\bf \S 3. Proofs\end{center}

\textit{Proof of Theorem }1. Let $\{B_n\}_{n=1}^{\infty}$ be an
exhaustion of the domain $f(D)$ by simply connected domains $B_n$,
where $\infty\in B_n$, $\overline{B}_n\subset
B_{n+1},\;n=1,2,\ldots$, and $\bigcup\limits_{n=1}^{\infty}
B_n=f(D)$. For any $n$ and sufficiently large $v>0$, the function
$u_n(w):=g_{B_n}(w,iv)-g_{B_n}(w,\infty)$ is harmonic in the
domain $B_n\setminus \{iv,\infty\}$, and we have
$u_n(w)\rightarrow +\infty$ as $w\rightarrow iv$ and
$u_n(w)\rightarrow -\infty$ as $w\rightarrow\infty$. In
particular, this implies that the sets $B_n^1:=\{w:u_n(w)>0\}$ and
$B_n^2:=\{w:u_n(w)<0\}$ are disjoint domains. The Green function
of the domain $B_n^1$ with a pole at the point $iv$ coincides with
the function $u_n(w)$ on $B_n^1$. Hence,
$$
\log\,r(B_n^1,iv)=\lim\limits_{w\rightarrow iv}(u_n(w)+\log
|w-iv|)=\log\,r(B_n,iv)-g_{B_n}(iv,\infty).
$$
Similarly, the function $-u_n(w)$ coincides on $B_n^2$ with the
Green function of this domain with a pole at the point $w=\infty$.
Therefore,
$$
\log\,r(B_n^2,\infty)=\lim\limits_{w\rightarrow\infty}(-u_n(w)
-\log|w|)=\log\,r(B_n,\infty)-g_{B_n}(\infty,iv).
$$
Adding these relations, we obtain
$$
\log[r(B_n,iv)r(B_n,\infty)]-2g_{B_n}(iv,\infty)=\log
[r(B_n^1,iv)r(B_n^2,\infty)]=
$$
$$
=\log[r(B_n^1\cup
B_n^2,iv)r(B_n^1\cup B_n^2,\infty)]+ 2g_{B_n^1\cup
B_n^2}(iv,\infty).\eqno(4)
$$
By Lemma 2, the last expression does not exceed the sum
$$
\log[r(S_v(B_n^1\cup B_n^2),iv)r(S_v(B_n^1\cup B_n^2),\infty)+2g_{
S_v(B_n^1\cup B_n^2)}(iv,\infty).
$$
Since $B_n^1\cap B_n^2=\emptyset$, it follows that for any
$\theta$ the ray $w=iv+\rho
e^{i\theta},\;0\leqslant\rho\leqslant\infty$, meets the set
$\overline{\mathbb{C}}_w\setminus(B_n^1\cup B_n^2)$. Hence, by the
definition of symmetrization with respect to a circle, the set
$R_v(\overline{\mathbb{C}}_w\setminus(B_n^1\cup B_n^2))$ contains
the circle $|w-iv|=v$. Therefore, the connected components
$\widetilde{B}_n^1$ and $\widetilde{B}_n^2$ of the set
$S_v(B_n^1\cup B_n^2)$ that contain the points $iv$ and
$w=\infty$, respectively, are disjoint, and we have
$$
r(S_v(B_n^1\cup B_n^2),iv)=r( \widetilde{B}_n^1,iv),\,\,\,\,\,\,\
r(S_v(B_n^1\cup B_n^2),\infty)=r( \widetilde{B}_n^2,\infty),
$$
$$
g_{S_v(B_n^1\cup B_n^2)}(iv,\infty)=g_{ \widetilde{B}_n^1\cup
\widetilde{B}_n^2}(iv,\infty)=0.
$$

We finally obtain
$$
\log[r(B_n,iv)r(B_n,\infty)]-2g_{B_n}(iv,\infty) \leqslant\log[r(
\widetilde{B}_n^1 ,iv)r(
\widetilde{B}_n^2,\infty)-2g_{\widetilde{B}_n^1\cup
\widetilde{B}_n^2}(iv,\infty)=$$
$$
=\log[r(S_vB_n,iv)r(S_vB_n,\infty)]-2g_{S_vB_n}(iv,\infty).\eqno(5)
$$
The last equation is established in the same way as (4) in view of
the symmetry of the domain $S_vB_n$ with respect to the circle
$|w-iv|=v$. Further, we shall use the following fact, which was
proved for the first time by Schiffer [4, \S3]: if domains $G_1$
and $G_2$ admit Green functions, $G_1\subset G_2$, and $\zeta$ and
$w$ are distinct points of $G_1$, then
$$
\log[r(G_1,\zeta)r(G_1,w)]-2g_{G_1}(\zeta,w)\leqslant
\log[r(G_2,\zeta)r(G_2,w)]-2g_{G_2}(\zeta,w).\eqno(6)
$$
By Lemma 1,
$$
S_vB_n\subset\overline{\mathbb{C}}_w\setminus(
\overline{\mathbb{C}}_w\setminus f(D))^*\subset f^*(D).
$$
Thus, by (5) and (6), the following inequality holds:
$$
\log[r(B_n,iv)r(B_n,\infty)]-2g_{B_n}(iv,\infty)\leqslant$$
$$ \leqslant
\log[r(f^*(D),iv)r(f^*(D),\infty)]-2g_{f^*(D)}(iv,\infty).\eqno(7)
$$
We denote by
$$
f_n(z)=a_1^nz+a_0^n+\frac{a_{-1}^n}{z}+\ldots
$$
a function which maps $D$ conformally and univalently onto the
domain $B_n$. Let
$$
h_n(w)=\frac{w}{a_1^n}-\frac{a_0^n}{a_1^n}-\frac{a_{-1}^n}{w}+\ldots
$$
be the expansion of the inverse map in a neighbourhood of the
point at infinity. Then the following equations hold for any
sufficiently large $v$:
$$
r(B_n,\infty)|a_1^n|=r(D,\infty)=1,\,\,\,\,\,r(B_n,iv)|h'_n(iv)|=r(D,h_n(iv))=|h_n(iv)|^2-1,
$$
$$
g_{B_n}(iv,\infty)=g_D(h_n(iv),\infty)=\log|h_n(iv)|.
$$
This gives
$$
r(B_n,iv)r(B_n,\infty)e^{-2g_{B_n}(iv,\infty)}=\frac{
|h_n(iv)|^2-1}{|a_1^n\;\;h'_n(iv)||h_n(iv)|^2}=$$
$$
=\frac{1-\left|\dfrac{iv}{a_1^n}+O(1)\right|^{-2}}{
\left|1-\dfrac{a_1^na_{-1}^n}{v^2}+o\left(\dfrac{1}{v^2}\right)\right|}
=1-(|a_1^n|^2-\mathrm{Re}\,a_1^na_{-1}^n)\frac{1}{v^2}+o\left(
\frac{1}{v^2}\right),\;v\rightarrow +\infty.
$$
Repeating the above manipulations with the function $f^*$ instead
of $f_n$, we see from the inequality (7) that
$$
|a_1^n|^2-\mathrm{Re}\,a_1^na_{-1}^n\geqslant
|a_1^*|^2-\mathrm{Re}\,a_1^*a_{-1}^*.
$$
Passing to the limit as $n\rightarrow\infty$, we obtain the inequality (2). This completes the
proof of the theorem.

\textit{Proof of Theorem }2. This follows from inequality (6),
where one must set $G_1=f^*(D),\;G_2=\widetilde{f}(D),\;\zeta=iv$
and $w=\infty$, and where $v$ is sufficiently large. By repeating
for $f^*$ and $\widetilde{f}$ the calculations for $f_n$ in the
last part of the proof of Theorem 1, we finally obtain the
inequality (3).

\textit{Proof of Corollary }1. We first consider the case in which
$w_0=0$ and $\varphi=\pi/2$. The function
$$
\widetilde{f}(z):=iR_f(0)h^{-1}(\lambda h(z))=i\lambda
R_f(0)z+\frac{iR_f(0)}{z}\left(\lambda-\frac{1}{\lambda}\right)+\ldots
$$
maps D conformally and univalently onto the complement of the set
$\widetilde{E}:=\{w:|w|\leqslant R_f(0)\}\cup
\{w:\mathrm{Re}\,w=0,\;|\mathrm{Im}\,w|\leqslant m_f(0,\pi/2)/2\}$
when
$\lambda=[m_f^2(0,\pi/2)/(4R_f^2(0))+1]/[m_f(0,\pi/2)/R_f(0)]^{-1}$.
Under the hypotheses of the corollary, we have $
\widetilde{E}\subset E^*$ ($E=\overline{\mathbb{C}}_w\setminus
f(D)$). Therefore, it follows from the inequalities (2) and (3)
that
$$
1-\mathrm{Re}\,a_{-1}\geqslant R_f^2(0)(2\lambda^2-1)=\frac{
m_f^4(0,\pi/2)+16R_f^4(0)}{8m_f^2(0,\pi/2)}.
$$
If $w_0$ and $\varphi$ are arbitrary, one must apply the above
conclusion to the function
$$
e^{i(\pi/2-\varphi)}[f(e^{i(\varphi-\pi/2)}z)-w_0]=z+\widetilde{a}_0-
e^{-2i\varphi}a_{-1}/z+\ldots\;.
$$
The condition for the equation can be verified immediately.

\textit{Proof of Corollary }2. We may assume that $\beta=-\gamma$.
The function
$$
F(\zeta)=c\int\limits_{0}^{\zeta}\sqrt{\frac{\zeta^2-k^2}{\zeta^2-1}}
d\zeta+\frac{i\alpha}{2}=c\zeta+c_0-\frac{c(1-k^2)}{2\zeta}+\ldots
$$
maps the upper half-plane $\mathrm{Im}\,\zeta>0$ conformally and
univalently onto the quadrangle $\{w:\mathrm{Im}\,w>0\}\setminus
\widetilde{E}$, where we now have
$\widetilde{E}=\{w:|\mathrm{Re}\,w|\leqslant\gamma,\;|\mathrm{Im}\,w|\leqslant\alpha/2\}$.
By the Riemann–Schwarz symmetry principle, the function
$$
\widetilde{f}(z):=F\left(\frac{1}{2}(z+\frac{1}{z})\right)=\frac{c}{2}z+
c_0-c\left(\frac{1}{2}-k^2\right)\frac{1}{z}+\ldots
$$
maps $D$ conformally and univalently onto the exterior of the
rectangle $\widetilde{E}$. Under the hypotheses of the corollary,
we have $\widetilde{E}\subset ( \overline{\mathbb{C}}_w\setminus
f(D))^*$. It remains to use the inequalities (2) and (3).

\begin{center}\bf Bibliography\end{center}

\begin{enumerate}
\item G. P\'{o}lya and G. Szeg\H{o}, \textit{Isoperimetric
inequalities in mathematical physics}, Ann. of Math. Stud., vol.
27, Princeton Univ. Press, Princeton, NJ 1951; Russian transl.,
Gosudarstv. Izdat. Fiz.-Mat. Lit., Moscow 1962.
\item G. F. Lawler, \textit{Conformally invariant processes in the plane}, Math. Surveys
Monogr., vol. 114, Amer. Math. Soc., Providence, RI 2005. 7
\item V.N. Dubinin, “Symmetrization in the geometric theory of functions
of a complex variable”, \textit{Uspekhi Mat. Nauk} \textbf{49}:1
(1994), 3–76; English transl., \textit{Russian Math. Surveys}
\textbf{49}:1 (1994), 1–79.
\item M.Schiffer, “Some new results in the theory of conformal mappings”,
Appendix to the book: R. Courant, \textit{Dirichlet’s Principle,
conformal mappings, and minimal surfaces}, Interscience, New York
1950; Russian transl., Inostr. Lit., Moscow 1953.
\item G. M. Goluzin, \textit{Geometric theory of functions of a complex variable}, Nauka, Moscow
1966; English transl., Transl. Math. Monogr., vol. 26, Amer. Math.
Soc., Providence, RI 1969.
\item W. K. Hayman, \textit{Multivalent functions}, Cambridge Univ. Press, Cambridge 1958; Russian transl.,
Inostr. Lit., Moscow 1960.
\item V. N. Dubinin, “Some properties
of the reduced inner modulus”, \textit{Sibirsk. Mat. Zh.}
\textbf{35}:4 (1994), 774–792; English transl., \textit{Siberian
Math. J.} \textbf{35}:4 (1994), 689–705.

\end{enumerate}

\end{document}